\documentclass[a4paper,10pt]{article}

\usepackage{amssymb}

\usepackage{graphicx}
\usepackage{graphics}

\def\yen{\hbox{iftdir\yoko\fi
\setbox0=\hbox{Y}Y\kern-.97\wd0\vbox{\hrule height.lex width.98\wd0
\kern.33ex\hrule height.lex width.98\wd0\kern.45ex}}}

\def\yen{{\setbox0=\hbox{Y}Y\kern-.97\wd0\vbox{hrule height.lex width.98%
\wd0\kern.33ex\hrule height.lex width.98\wd0\kern.45ex}}}
    
\def\F{$$\mathrm{ Figure }} 
\def\g{\gamma}
\def\g{\gamma}
\def\ve{\varepsilon}  
\def\vp{\varphi}
\def\si{\mathrm{sin}}
\def\co{\mathrm{cos}}
\def\i{\longrightarrow}
\def\e{\hookrightarrow}
\def\l{\longrightarrow} 
\def\ttt{\longmapsto}

\def\g{\gamma}
\def\ve{\varepsilon}
\def\vp{\varphi} 
 \def\si{\mathrm{sin}}
\def\co{\mathrm{cos}}
\def\i{\longrightarrow}
\def\e{\hookrightarrow}
\def\l{\longrightarrow} 
\def\ttt{\longmapsto}
\def\pa{ pair of $n$-knots}
\def\f{\flushpar }
\def\nl{\newline }
\def\np{\newpage }
\def\x{\times }
\def\p{\bf Proof of  }
\def\te{^t \hskip-1mm }

\begin{document}

\title{
The intersection of spheres in a sphere and \\ 
a new geometric meaning of the Arf invariant
}
\author{
 Eiji Ogasa\\
 ogasa@ms.u-tokyo.ac.jp\\
Department of Mathematical Sciences, 
 University of Tokyo\\ 
 Komaba, Tokyo 153, JAPAN  
 \thanks{ 
 Keywords: the Arf invariant, spin cobordism group, pass-moves, 
realizable pair of links, surgery of submanifolds \nl
This research was partially supported by Research Fellowships
of the Promotion of Science for Young Scientists.}
}
\date{}
\maketitle

\noindent{\bf Abstract}

Let $S^3_i$ be a 3-sphere embedded in the 5-sphere $S^5$ ($i=1,2$). 
Let $S^3_1$ and $S^3_2$ intersect transversely. 
Then the intersection $C=S^3_1\cap S^3_2$  
is a disjoint collection of circles.  
Thus we obtain a pair of 1-links,   $C$ 
in $S^3_i$ ($i=1,2$), 
and a pair of  3-knots, $S^3_i$ in $S^5$ ($i=1,2$). 
Conversely let $(L_1,L_2)$ be a pair of 1-links and $(X_1,X_2)$ be a pair 
of 3-knots. 
It is natural to ask whether 
the pair of 1-links $(L_1,L_2)$    
is obtained as the intersection of 
the 3-knots $X_1$ and $X_2$  
as above. 
We give a complete answer to this question. 
Our answer gives a new geometric meaning of the Arf invariant 
of 1-links. 

Let $f:S^3 \longrightarrow S^5$ be a smooth transverse immersion.  
Then the self-intersection  $C$ consists of double points.  
Suppose that $C$ is a single circle in $S^5$. 
Then $f^{-1}(C)$  in $S^3$ is a 1-knot or a 2-component 1-link.   
There is a similar realization problem. 
We give a complete answer to this question.

\np

\section{  Introduction and Main results } 

Let $S^3_i$ be a 3-sphere embedded in the 5-sphere $S^5$ ($i=1,2$). 
Let $S^3_1$ and $S^3_2$ intersect transversely. 
Then the intersection $C=S^3_1\cap S^3_2$  
is a disjoint collection of circles. 
Then $C$ in $S^3_i$ is a 1-link ($i=1,2$).
Note that the orientation of $C$ is induced by 
that of $S^3_1$ , 
that of $S^3_2$ and that of $S^5$.  
Thus we obtain a pair of 1-links,   $C$ 
in $S^3_i$ ($i=1,2$), 
and a pair of  3-knots, $S^3_i$ in $S^5$ ($i=1,2$). 

Conversely let $(L_1,L_2)$ be a pair of 1-links and $(X_1,X_2)$ be a pair 
of 3-knots. It is natural to ask whether 
 the pair of 1-links $(L_1,L_2)$ 
is obtained as the 
intersection of 
the 3-knots $X_1$ and $X_2$ 
as above. 
We give a complete answer to this question. (Theorem 1.1.)

To state our results we need some definitions.    
 
An  {\it (oriented) (ordered) $m$-component n-(dimensional) link}
 is a smooth, oriented submanifold $L=\{K_1,...,K_m\}$ of $S^{n+2}$, 
which is the ordered disjoint union of $m$ manifolds, each PL homeomorphic 
to the $n$-sphere.  
If $m=1$, then $L$ is called a {\it knot}.  
(See [1], [2], [8], [9]. )
 
We say that n-links $L_1$ and $L_2$ are {\it equivalent} 
if there exists an orientation preserving diffeomorphism 
$f:$ $S^{n+2}$ $\rightarrow$ $S^{n+2}$ 
such that $f(L_1)$=$L_2$  and 
$f\vert_{L_1}:$ $L_1$ $\rightarrow$ $L_2$ is 
an orientation preserving diffeomorphism.

{\bf Definition  }  
$(L_1, L_2, X_1, X_2)$ is called a 
{\it 4-tuple of links } if the following conditions (1), (2) and (3) hold.

(1)  
  $L_i=(K_{i1},...,K_{im_i})$  is an oriented ordered $m_i$-component 
1-dimensional link  $(i=1,2).$  

(2)  
   $m_1=m_2.$

(3)  
    $X_i$ is a 3-knot.

{\bf Definition   } 
A 4-tuple of links $(L_1, L_2, X_1, X_2)$ is said to be  {\it realizable } if 
there exists a smooth transverse 
immersion  $f:S^3_1\coprod S^3_2 \longrightarrow S^5$   
with the  following properties.  
We assume that the orientations of $S^3_1, S^3_2$ 
and $S^5$ are given.

(1)  
 $f\vert S^3_i$ is a smooth embedding. 
  $f(S^3_i)$ in $S^5$ is equivalent to the  3-knot $X_i (i=1,2).$  

(2)
 For  $C=f(S^3_1)\cap f(S^3_2)$, 
the inverse image  $f^{-1}(C)$ in $S^3_i$ is equivalent to the 1-link 
$L_i (i=1,2).$   
Here, the orientation of $C$ is induced 
naturally from the preferred orientations of $S^3_1, S^3_2,$ and $S^5$,   
 and an arbitrary order is given to the components of $C$.

The following theorem characterizes the realizable 4-tuples of links.

{\bf Theorem 1.1}    {\it 
 A 4-tuple of links $(L_1, L_2, X_1, X_2)$ is realizable  if and only if 
 $(L_1, L_2, X_1, X_2)$  satisfies one of the following conditions (1) and (2).

(1)
Both $L_1$ and $L_2$ are proper links, and 

\vskip3mm
\hskip2cm $\mathrm{Arf}(L_1) = \mathrm{Arf}(L_2).$
\vskip3mm

(2)
 Neither $L_1$ nor $L_2$ is a proper link, and 

\vskip3mm
\hskip2cm
$\mathrm{lk}(K_{1j}, L_1-K_{1j})\equiv
  \mathrm{lk}(K_{2j}, L_2-K_{2j})  \mathrm{mod 2}$ 
  for  all $j$. 
}
\vskip3mm

In the case where $L_i$ is a 1-component link, 
that is, $L_i$ is a knot $K_i$, we have the following corollary. 

{\bf Corollary 1.2}   {\it 
 For 1-knots  $K_1$ and $K_2$, a 4-tuple of links  
 $(K_1, K_2, X_1, X_2)$ is realizable    
if and only if

\vskip3mm
\hskip2cm $\mathrm{Arf}(K_1) = \mathrm{Arf}(K_2).$
\vskip3mm
}

{\bf Note.}   
Theorem 1.1 and Corollary 1.2 give 
a new geometric meaning of the Arf invariant. 

{\bf Note 1.2.1.}     
The  problem (26) in [3]
says: Investigate ordinary sense slice 1-links,  
where {\it ordinary sense slice 1-links} 
are 1-links which are obtained as follows: 
Let $S^2$ be in 
$ R^4= R^3\times R.$ 
Then 
$S^2\cap [ R^3\times\{0\}]$ in $ R^3\times\{0\}$ is a 1-link. 
By using Theorem 1.1, the author 
gives  an  answer to this  problem.    
The answer is: 
for every ordinary sense slice 1-link     
we can define the Arf invariant and it is zero (see [18]).

Let $f:S^3 \longrightarrow S^5$ be a smooth transverse immersion.  
Then the self-intersection  $C$ consists of double points.  
Suppose that $C$ is a single circle  in $S^5$. 
Then the $f^{-1}(C)$  in $S^3$ is a 1-knot or a 2-component 1-link.   
There is a similar realization problem. 
We consider which 1-knots 
(resp. 2-component 1-links ) we obtain as above. 
We give complete answers.

{\bf Theorem 1.3}    {\it 
 Let $f:S^3 \longrightarrow S^5$ be a smooth transverse immersion.  
Then the self-intersection  $C$ consists of double points.  
Suppose that $C$ is a single circle  in $S^5$. 

(1)    Any 2-component 1-link is realizable 
as $f^{-1}(C)$  in $S^3$ for an immersion $f$.   

(2)    Any 1-knot is realizable 
as $f^{-1}(C)$  in $S^3$ for an immersion $f$.   
}

{\bf Remark. }  
Suppose that $K_1$ is the trivial 1-knot, 
$K_2$ is the trefoil 1-knot, and $X_1$ and $X_2$ are 3-knots.  
Suppose that $L=(L_1, L_2)$ is a split 1-link such that 
$L_1$ is the trivial 1-knot and $L_2$ is the trefoil knot.   
 Then, by Corollary 1.2,  a 4-tuple of links $(K_1, K_2, X_1, X_2)$   
is not realizable.  
But, by Theorem 1.3 (1),  the two component split link $L$ 
 is realizable as in Theorem 1.3 (1).

In [16],[17],[18],[19] 
the author discussed some topics which are related to this paper. 
In [19] he discussed the intersection of three 4-spheres in a 6-sphere. 
In [17] he discussed the intersection of two $(n+2)$-spheres in an $(n+4)$-sphere. 
In [16] he discussed the following: 
Let $L=(K_1, K_2)$ be a 2-link in $S^4=\partial B^5$. 
Take a slice disc $D^3_i$ in $B^5$  for each component $K_i$. 
He discussed the intersection of two slice discs $D^3_1$ and $D^3_2$ in the 5-ball $B^5$.  
In [18] he applied Theorem 1.1 to Fox's problem as we state in Note 1.2.1.

{\bf Problem 1.4. } 
Suppose each of $S^3_1, S^3_2$, and  $S^3_3$ is a 3-sphere embedded in $S^5$.  
Suppose $S^3_i$ and $S^3_j$ intersect transversely. 
Suppose each of $S^3_1\cap S^3_2$, $S^3_2\cap S^3_3$, and $S^3_3\cap S^3_1$ 
is a single circle. 
Then we have a triple of 1-links, 
$L_i$=($S^3_i\cap S^3_j$, $S^3_i\cap S^3_k$) in $S^3_i$,  
where $(i,j,k)=(1,2,3), (2,3,1), (3,1,2)$. 

Which triple of 1-links do we obtain like this?

Do we characterize such triple by the Arf invariants, the linking numbers, and 
the Saito-Sato-Levine invariants? 
(See [20] for the definition of the Saito-Sato-Levine invariant.)

In [19] the author discussed a higher dimensional version of Problem 1.4.

This paper is organized as follows. 
In \S 2 
we review  
spin cobordism and the Arf invariant. 
In \S 3
 we discuss 
a necessary condition 
for the realization of 4-tuple of links. 
We find the obstruction for the realization in the spin cobordism group
$\Omega^{\mathrm{spin}}_{\star}.$ 
 In \S 4 
 we discuss a sufficient condition for the realization of 4-tuple of links. 
We carry out surgeries of submanifolds to carry out an (un)knotting  operation. 
Theorem 1.1 is deduced from \S 3 and \S 4. 
In \S 5 
we prove Theorem 1.3. 
In \S 6 
we give a problem.


\section{ 
Spin cobordism and the Arf invariant}

In this section we review some results on 
the Arf invariant and spin cobordism. 
See [5] and [6] 
for the Arf invariant. 
See [10] 
for spin structures and spin cobordism.

We suppose that, when we say $M$ is  a spin manifold,  
$M$ is oriented.

Recall that a {\it proper link} is an $m$-component 1-link    $L=\{K_1,...,K_m\}$ 
such that 
$ \mathrm{lk}(K_j, L-K_j)$ $=\sum_{1\leqq i \leqq m, i\not= j}$ 
$ \mathrm{lk}(K_j, K_i)$ 
is an even number for each $K_j$.

Let $L=(K_1,...,K_m)$ be an $m$-component 1-link. 
Let $F$ be a Seifert surface for $L$.  
We induce a spin structure $\sigma$ on $F$ from the unique one on $S^3$.  
We induce a spin structure $\sigma_i$ on $K_i$ from $\sigma$ on $F$.  
Then we have: 

{\bf Proposition 2.1}   {\it 
Under the above condition, 
for each $i$, 

mod 2  $\mathrm{lk}$ $(K_i, L-K_i)\quad
=[(K_i, \sigma_i)]\in\Omega_1^{\mathrm{spin}}$.   

In particular, $L$ is a proper link if and only if each $[(K_i, \sigma_i)]=0$. 
}

Suppose that $L$ is a proper link. Take $(F, \sigma)$ as above. 
Let $\hat{F}$ be the closed  surface obtained 
from $F$ by attaching disks to the boundaries. 
Let $\hat{\sigma}$ be the unique extension of 
$\sigma$ over $\hat{F}$. 
Then we have:

{\bf Proposition 2.2}   {\it 
  Under the above condition, 
 $\mathrm{Arf}$ $(L)$= $[(\hat{F}, \hat{\sigma})]$ $\in\Omega_2^{\mathrm{spin}}$.  }

Although they may be folklore, 
the author gives a proof of Proposition 2.1 and 
that of Proposition 2.2 in the appendix.

\section{ 
A necessary condition for the realization of 4-tuple of links  
}

In this section we discuss a necessary condition 
for the realization of a 4-tuple of links. That is, 
we prove the  following two propositions. 
  
{\bf Proposition 3.1}   {\it 
If $(L_1, L_2, X_1, X_2)$ is realizable then 

$\hskip36pt
\ \mathrm{lk}(K_{1j}, L_1-K_{1j})
\equiv
\ \mathrm{lk}(K_{2j}, L_2-K_{2j}) \quad \mathrm{mod}\hskip2pt2
$\quad
for all $j$. 
\newline
In particular,  
$L_1$ is a proper link 
if and only if 
$L_2$ is a proper link.  
}

{\bf Proposition 3.2}   {\it 
Let $L_1$ and $L_2$ be proper links. 
If $(L_1, L_2, X_1, X_2)$
is realizable then 
$$\mathrm{Arf}(L_1) = \mathrm{Arf}(L_2).$$ 
}

In order to prove them, we prepare a lemma.

Let $M_1$ and  $M_2$ be codimension one submanifolds of 
an $n$-dimensional compact spin manifold $N$. 
Suppose that $M_1$ and $M_2$ are compact oriented manifolds. 
Suppose that $M_1$, $M_2$, and $N$ may have the boundary and the corner.  
Let $M_1$ and  $M_2$  intersect transversely.  
Suppose $M_i$ may be embedded in the boundary (resp. the corner ) of $N$. 

We induce a spin structure $\sigma_i$ on $M_i$ from $N$. 
We induce a spin structure $\xi_i$ on $M_1\cap M_2$ from  $\sigma_i$ on $M_i$.  
($i=1,2$). 
Then it is easy to prove:

{\bf Lemma }   {\it 
$\xi_1$ and $\xi_2$ are same.  
} 

The spin structure $\xi_1$=$\xi_2$ on   $M_1\cap M_2$ is called 
 {\it the unique spin structure induced by $M_1$, $M_2$ and $N$}. 
 

\vskip1cm
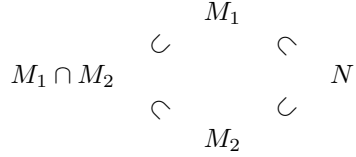
\begin{tabular} {llllc}
&&$M_1$&&\\
&
\rotatebox[origin=c]{45}{$\subset$}
&&
\rotatebox[origin=c]{-45}{$\subset$}
&\\
$M_1\cap M_2$&&&&$N$\\
&
\rotatebox[origin=c]{-45}{$\subset$}
&&
\rotatebox[origin=c]{45}{$\subset$}
&\\
&&$M_2$&&\\

\end{tabular}
\vskip1cm


{\bf Proof of Proposition 3.1.}
Let $f:S^3_1\amalg S^3_2\i S^5$ be an immersion to realize 
the 4-tuple of links $(L_1, L_2, X_1, X_2)$.  
Let $C=\amalg C_i$ denote $f(S^3_1)\cap f(S^3_2)$.

We abbreviate $f(S^3_i)$ to  $ S^3_i.$   
Let $V_i$ be a Seifert hypersurface for $X_i$, i.e.,  

$$S^3_i=\partial V_i\subset V_i\subset S^5.$$

We  make $V_1$ and $V_2$ intersect transversely.  
We induce a spin structure $v_i$ on $V_i$ from the unique one on $S^5$.

Put $W=V_1\cap V_2$. 
We have: 


\vskip1cm
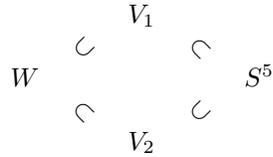
\begin{tabular} {llllc}
&&$V_1$&&\\
&
\rotatebox[origin=c]{45}{$\subset$}
&&
\rotatebox[origin=c]{-45}{$\subset$}
&\\
$W$&&&&$S^5$\\
&
\rotatebox[origin=c]{-45}{$\subset$}
&&
\rotatebox[origin=c]{45}{$\subset$}
&\\
&&$V_2$&&\\
\end{tabular}
\vskip1cm


We give $W$ the unique spin structure $w$ induced by $V_1$, $V_2$ and $W$.

Here, we see that 
$\partial W=(\partial V_1 \cap V_2)\cup (V_1 \cap \partial V_2).$  
Put  $F_1$=$\partial V_1 \cap V_2$.  
Put  $F_2$=$V_1 \cap \partial V_2$. 
Then $F_1$ ( resp. $F_2$ ) is a Seifert surface for $L_1$ ( resp. $L_2$ ). 

We induce a spin structure on $\partial V_i=S^3_i$ from $v_i$ on $V_i$. 
Note that it is the unique one on $S^3_i$.  

We have:


\vskip1cm
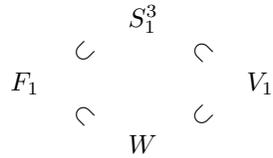
\begin{tabular} {llllc}
&&$S^3_1$&&\\
&
\rotatebox[origin=c]{45}{$\subset$}
&&
\rotatebox[origin=c]{-45}{$\subset$}
&\\
$F_1$&&&&$V_1$\\
&
\rotatebox[origin=c]{-45}{$\subset$}
&&
\rotatebox[origin=c]{45}{$\subset$}
&\\
&&$W$&&\\
\end{tabular}
\vskip1cm


We give  $F_1$ the unique spin structure $\sigma_1$ induced by $S^3_1$, $W$ and $V_1$.

We have:


\vskip1cm
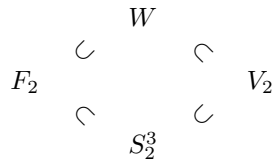
\begin{tabular} {llllc}
&&$W$&&\\
&
\rotatebox[origin=c]{45}{$\subset$}
&&
\rotatebox[origin=c]{-45}{$\subset$}
&\\
$F_2$&&&&$V_2$\\
&
\rotatebox[origin=c]{-45}{$\subset$}
&&
\rotatebox[origin=c]{45}{$\subset$}
&\\
&&$S^3_2$&&\\

\end{tabular}
\vskip1cm


We give  $F_2$ the unique spin structure $\sigma_2$ induced by $S^3_2$, $W$ and $V_2$.

We have:


\vskip1cm
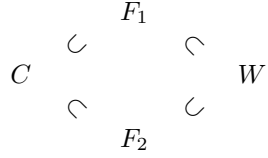
\begin{tabular} {llllc}
&&$F_1$&&\\
&
\rotatebox[origin=c]{45}{$\subset$}
&&
\rotatebox[origin=c]{-45}{$\subset$}
&\\
$C$&&&&$W$\\
&
\rotatebox[origin=c]{-45}{$\subset$}
&&
\rotatebox[origin=c]{45}{$\subset$}
&\\
&&$F_2$&&\\
\end{tabular}
\vskip1cm


We give $C_j$ the unique spin structure  $\tau_j$ induced 
by $F_1$, $F_2$ and $W$.

Then 
 mod2 lk$(K_{1j}, L_1-K_{1j})$=mod2 lk$(K_{2j}, L_2-K_{2j})$  
=$[(C_j, \tau_j)]$ $\in\Omega_1^{\mathrm{spin}}$  
for all $j$. 
This completes the proof of Proposition 3.1. 

We confirm that we have:


\vskip1cm
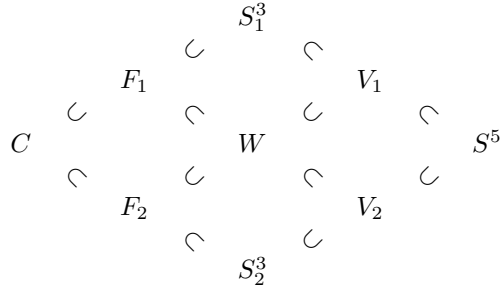
\begin{tabular} {llllllllc}
&&&&$S^3_1$&&&&\\
&&&\rotatebox[origin=c]{45}{$\subset$}
&&
\rotatebox[origin=c]{-45}{$\subset$}
&&&\\
&&$F_1$&&&&$V_1$&&\\
&
\rotatebox[origin=c]{45}{$\subset$}
&&
\rotatebox[origin=c]{-45}{$\subset$}
&&
\rotatebox[origin=c]{45}{$\subset$}
&&
\rotatebox[origin=c]{-45}{$\subset$}
&\\
$C$&&&&$W$&&&&$S^5$\\
&
\rotatebox[origin=c]{-45}{$\subset$}
&&
\rotatebox[origin=c]{45}{$\subset$}
&&
\rotatebox[origin=c]{-45}{$\subset$}
&&
\rotatebox[origin=c]{45}{$\subset$}
&\\
&&$F_2$&&&&$V_2$&&\\
&&&\rotatebox[origin=c]{-45}{$\subset$}
&&
\rotatebox[origin=c]{45}{$\subset$}
&&&\\
&&&&$S^3_2$&&&&\\
\end{tabular}
\vskip1cm


{\bf Proof of Proposition 3.2.}
We can make the corner of $W$ smooth. 
Note 
$\partial (W, w)$ =
$({F_1},  {\sigma_1})\cup ({F_2},  {\sigma_2}$. 
Take  $\hat{F_i}$ and $\hat{\sigma_i}$ as above. 
Then we have: 

$\mathrm{Arf}(L_1) + \mathrm{Arf}(-L_2)=
\mathrm{Arf}(L_1) + \mathrm{Arf}(L_2)=
[(\hat{F_1},  \hat{\sigma_1})]
+[(\hat{F_2},  \hat{\sigma_2})]
=[\partial (W, w)] 
=0\in \Omega^{\mathrm{spin}}_2\cong Z_2. $

Therefore  
$\mathrm{Arf}(L_1) = \mathrm{Arf}(L_2).$

\section{ 
 A sufficient condition for the realization of 4-tuple of links  
}

In this section we discuss 
a sufficient condition 
for the realization of a 4-tuple of links.  
That is, 
we prove the following proposition. 

{\bf Proposition 4.1}  {\it 
A  4-tuple of links $(L_1, L_2, X_1, X_2)$ is realizable  
if $(L_1, L_2, X_1, X_2)$  satisfies one of 
the  conditions (1) and (2) of Theorem 1.1.  
}

It is easy to prove that Proposition 4.1 is equivalent to the following Proposition 4.2.

{\bf Proposition 4.2 }  {\it 
Let $X_1$ and $X_2$ be the trivial 3-knots. 
Let $L_1$ and $L_2$ be 1-links.  
Suppose that $L_1$ and $L_2$  satisfies 
one of the  conditions (1) and (2) of Theorem 1.1.  
Then the 4-tuple of links $(L_1, L_2, X_1, X_2)$ 
is realizable.   
}

We prove 

{\bf  Lemma 4.3 } {\it      
Let $X_1$ and $X_2$ be the trivial 3-knots. 
Let $L$ be a 1-link.   
Then the 4-tuple of links $(L, L, X_1, X_2)$ is realizable.   
}

{\bf Proof of Lemma 4.3. }
 Let $f:S^3_1\coprod S^3_2 \e S^5$ 
be an embedding  such that  
$f (S^3_1\coprod S^3_2)$ in $S^5$ 
is equivalent to the  trivial 3-link. 
We  take a chart $(U, \phi)$  of $S^5$ 
with the following properties (1) and (2). 

(1)
 $\phi:U\cong  R^5 = \{(x,y,z,u,v)\vert x,y,z,u,v\in   R\}$
=$ R^3\x  R_u\x R_v$.

(2)
$U\cap f(S^3_1)=\{(x,y,z,u,v)\vert u=0, v=0\}= R^3_1$ 

$U\cap f(S^3_2)=\{(x,y,z,u,v)\vert u=1, v=0\}= R^3_2$




$$\mathrm{Figure 1}$$

Obviously Lemma 4.3 follows from Lemma 4.4.

{\bf  Lemma 4.4 }     {\it 
There exists an immersion  $g: S^3_1\coprod S^3_2 \longrightarrow S^5$ 
with the following conditions.

(1)
$g\vert_{ S^3_2} =f\vert_{ S^3_2}$.  

(2)
$g\vert_{ S^3_1}$ is isotopic to $f\vert_{ S^3_1}$. 

(3)
$g\vert_{ S^3_1-g^{-1}(U)}$ = $f\vert_{ S^3_1-g^{-1}(U)}$.   
Hence $g(S^3_1)\cap g(S^3_2)$ $\subset U$.

(4)
$g(S^3_1)\cap g(S^3_2)$ in 
$g(S^3_1)$ and that in $g(S^3_2)$ are both equivalent to the 1-link $L$.  

(5)
$g(S^3_1)\cap\{(x,y,z,u,v)\vert u=0, v\in R\}$= 
 $g(S^3_1)\cap\{(x,y,z,u,v)\vert u=0, v=0\}$.  
}

{\bf Proof of Lemma 4.4. }
We modify the embedding $f$ to construct an immersion  $g$.

In 
$ R^3_1$ 
we take the 1-link $L$ and a Seifert surface $F$ for $L$.     
Let $N(F)$

=$F\times\{t\vert -1\leqq t\leqq1\}$ 
be a tubular neighborhood of $F$ in 
$ R^3_1$.     

We define a subset $E$ of 
$N(F)\times  R_u \times  R_v $
$=\{(p,t,u,v)\vert p\in F, -1\leqq t\leqq 1, u\in R, v\in R \}$
so that  

$E=\{(p,t,u,v)\vert p\in F, \quad$ 
$ 0\leqq u\leqq \frac{\pi}{2},\quad 
t=k\cdot \mathrm{cos} u, \quad
v=k\cdot \mathrm{sin} u, \quad$ 
$-1\leqq k\leqq 1\}$.

Put  $P$=$\overline{\partial E-(\partial E\cap f(S^3_1))}$.
Put  $Q$=$\overline{f(S^3_1)-(\partial E\cap f(S^3_1))}$.
Note that $\partial P$=$\partial Q$=$\partial N(F)$.  
Put $\Sigma=$$P\cup Q$.
Then, by the construction, $\Sigma$ is a 3-sphere embedded in $S^5$ and 
is the trivial 3-knot.  

\vskip3mm
{\bf Note.} 
In $U$ the following hold. 

(1)   $g(S^3_1)\cap\{(x,y,z,u,v)\vert u=0, v\in R\}$ 
        $=g(S^3_1)\cap\{(x,y,z,u,v)\vert u=0, v=0\}$ 
        $=\overline{g(S^3_1)-N(F)}$

(2)   $g(S^3_1)\cap\{(x,y,z,u,v)\vert 0<u\leqq\frac{\pi}{2}, v\in R\}$ 
      is Int $N(F)$.

(3) Let $0<u'<\frac{\pi}{2}$.  
 
     $g(S^3_1)\cap\{(x,y,z,u,v)\vert u=u', v>0\}$ 
         is diffeomorphic to Int $F$. 
     
(4)  Let $0<u'<\frac{\pi}{2}$.  

     $g(S^3_1)\cap\{(x,y,z,u,v)\vert u=u', v<0\}$ 
     is diffeomorphic to Int $F$.

(5) Let $0<u'<\frac{\pi}{2}$. 

  $g(S^3_1)\cap\{(x,y,z,u,v)\vert u=u', v=0\}$ 
    is diffeomorphic to $L$.

(6)  Let $0<u'<\frac{\pi}{2}$.  

$g(S^3_1)\cap\{(x,y,z,u,v)\vert u=u', v\in R\}$ 
is diffeomorphic to $\partial(N(F))$. 

Let $F_i$ be diffeomorphic to $F$ ($i=1,2$). 
Recall $F$ is a compact oriented surface with boundary. 
We identify $\partial F_1$ with $\partial F_2$ 
to obtain $F_0=F_1\cup F_2$. 
Note $F_0$ is diffeomorphic to $\partial(N(F))$.

(7)    $g(S^3_1)\cap\{(x,y,z,u,v)\vert u=\frac{\pi}{2}, v\in R\}$ 
      is diffeomorphic to $N(F)$.

(8)    $g(S^3_1)\cap\{(x,y,z,u,v)\vert u=\frac{\pi}{2}, v=0\}$ 
      is diffeomorphic to $F$. 

(9)    $g(S^3_1)\cap\{(x,y,z,u,v)\vert 0<u\leqq\frac{\pi}{2}, v=0\}$ 
      is diffeomorphic to Int $F$.

\vskip3mm

In Figure 2, 3, 4,   
$\Sigma\cap U$ and $g(S^3_2)\cap U$  are drawn. 
There, we replace 
$ R^3\x R_u\x R_v$ with 
$ R^2\x R_u\x R_v$.  


$$\mathrm{Figure 2}$$

$$\mathrm{Figure 3}$$

$$\mathrm{Figure 4}$$

By the construction,  $\Sigma\cap$ $g(S^3_2)$  in  $\Sigma$ and 
that in $g(S^3_2)$  are both equivalent to $L$. 
Define $g\vert_{ S^3_1}$ so that $g(S^3_1)$=$\Sigma$.

This completes the proof of Lemma 4.4 and therefore Lemma 4.3.

{\bf Note.} 
 Lemma 4.3 gives an alternative proof of 
the results of [4]
and Lemma 1 of [15].

In order to prove Proposition 4.2, we review the pass-moves. 
See [5] and [6] for detail.

{\bf Definition } (See [5][6].)
Two 1-links are {\it pass-move equivalent}
 if one is obtained from the other 
by a sequence of pass-moves. 
See 
Figure 5 
for an illustrations of the pass-move. 


$$\mathrm{Figure 5}$$

In $B^3_a$ there are four arcs. 
Each of four arcs may belong to different components of the 1-link. 
We do not assume the four arcs belong to one component of the 1-link.

The following propositions are essentially proved in [5]. 
A proof is written in the appendix of [14]. 

{\bf Proposition 4.5 } (See [5][14].)  
   {\it 
Let $L_1$ and $L_2$ be 1-links. 
Then  $L_1$ and $L_2$ are pass-move equivalent if and only if 
$L_1$ and $L_2$
 satisfy 
one of the  conditions (1) and (2) of Theorem 1.1. 
}

{\bf Proposition 4.6 }   (See [5][14].)
 {\it 
Let $L_1$ and $L_2$ be 1-links. 
Suppose that $L_1$ is pass-move equivalent to $L_2$. 
Let $F$ be an oriented Seifert surface for $L_2$ 
such that the genus of $F$ is not less than  the genus of $L_1$. 
Then there exists a disjoint union of 3-balls $B^3$ 
such that  $B^3\cap F$ is  as  in 
Figure 6 
and the pass-moves in all $B^3$ change $L_2$ to $L_1$ 
}


$$\mathrm{Figure 6}$$

Let $L_1$ and $L_2$ be 1-links. 
Suppose $L_1$ and $L_2$ satisfy one of 
the conditions (1) and (2) of Theorem 1.1.  
Then by Proposition 4.5  
$L_1$ is obtained from $L_2$ by a sequence of pass-moves. 
We choose a Seifert surface $F$ for $L_2$ and 
the disjoint union of 3-balls $B^3$ in $S^3$ 
as in Proposition 4.6.    

Take a 1-link ($Y_{1}$, $Y_{2}$) in each $B^3$ as in 
Figure 7. 
By considering the Kirby moves of framed links (in [10]), 
it is easy to prove:

 
 $$\mathrm{Figure 7}$$

{\bf Lemma 4.7}  
 {\it 
We carry out 0-framed surgeries along all $Y_{1}$ and all $Y_{2}$.    
After these surgeries,  we have the following.  
(1)
$S^3$ becomes  a  3-sphere  again. 
(2)
each $B^3$ becomes a 3-ball again,  
(3) $L_2$ in the old sphere $S^3$ changes to $L_1$ in the new 3-sphere. 
} 

We go back to the proof of Proposition 4.2. 

As in Lemma 4.4, 
take an immersion $g: S^3_1\coprod S^3_2 \longrightarrow S^5$ 
to realize $(L_2, L_2, X_1, X_2)$, where $X_i$ is the trivial 3-knot.  

Obviously Proposition 4.2 follows from Lemma 4.8.

 {\bf Lemma 4.8 } 
 {\it 
There exists an immersion  $h: S^3_1\coprod S^3_2 \longrightarrow S^5$ 
with  the following conditions.

(1)
$h\vert_{ S^3_2} =g\vert_{ S^3_2}$.  

(2)
$h( S^3_1 ) $ in $S^5$ is equivalent to the trivial 3-knot.  

(3)
$h\vert_{ S^3_1-g^{-1}(U)}$ = $g\vert_{ S^3_1-g^{-1}(U)}$. 
Hence $h(S^3_1)\cap h(S^3_2)$ $\subset U$.

(4)
$h(S^3_1)\cap h(S^3_2)$ in $h(S^3_1)$ is equivalent to $L_1$. 
$h(S^3_1)\cap h(S^3_2)$ in $h(S^3_2)$ is equivalent to $L_2$. 

(5)
$g(S^3_1)\cap\{u\neq0, v\in R\}$=$h(S^3_1)\cap\{u\neq0, v\in R\}$  
}

{\bf Proof of Lemma 4.8. } 
We  modify the immersion $g$ to define an immersion $h$ as follows. 

As in the proof of Lemma 4.4, 
we take $L_2$ and a Seifert surface $F$ for $L_2$ in $g(S^3_1)$. 
Suppose 
$L_2\subset R^3_1$ and $F\subset U$. 
Suppose the genus of $F$ satisfies the condition in Proposition 4.6. 
Take 3-balls $B^3$  in $g(S^3_1)$  as in Proposition 4.6.   
Take ($Y_1$, $Y_2$) in $B^3$ as in Lemma 4.7.  
See Int $N(F)$ in $g(S^3_1)$. 
Int $N(F)$ is not in 
$ R^3_1$ 
and is in $\{u>0\}$. 
But by Proposition 4.6, 
we can suppose $Y_1$ and $Y_2$ are in 
$ R^3_1$.  

Take  an embedded 2-disc $h^{'2}_{1}$  in 
$\{(x,y,z,u,v)\vert u=0, v\geqq 0 \}$ 
so that 
$h^{'2}_{1}$  meets  $\{(x,y,z,u,v)\vert u=0, v=0 \}$ 
at $Y_1$ transversely. 
Suppose that $h^{'2}_{1}$ is embedded trivially.

Take  an embedded 2-disc $h^{'2}_{2}$  in 
$\{(x,y,z,u,v)\vert u=0, v\leqq 0 \}$ 
so that 
$h^{'2}_{2}$  meets  $\{(x,y,z,u,v)\vert u=0, v=0 \}$ 
at $Y_2$ transversely. 
Suppose that $h^{'2}_{2}$ is embedded trivially. 

See Figure 8.


$$\mathrm{Figure 8}$$

Let $h^2_{1}$ be a tubular neighborhood of $h^{'2}_{1}$ 
in $\{(x,y,z,u,v)\vert u=0, v\geqq 0 \}$.

Let $h^2_{2}$ be a tubular neighborhood of $h^{'2}_{2}$ 
in $\{(x,y,z,u,v)\vert u=0, v\leqq 0 \}$.  

Then $h^2_{1}$ and $h^2_{2}$ are diffeomorphic to the 4-ball.    
Furthermore we can regard $h^2_i$  
as a 4-dimensional 2-handle  attached to  
$g(S^3_1)$ along $Y_i$  with 0-framing.

Put $R=g(S^3_1)$ $-\{g(S^3_1)\cap\partial h^2_1\}$ 
$-\{g(S^3_1)\cap\partial h^2_2\}$.
Put $S_1=$ $\partial h^2_1$ $-\{g(S^3_1)\cap\partial h^2_1\}$.  
Put $S_2=$ $\partial h^2_2$ $-\{g(S^3_1)\cap\partial h^2_2\}$.  
Here, we consider all $h^2_i$.  
Put $\Lambda$=$R\cup S_1\cup S_2$.

Then we can regard  $\Lambda$ as the result of 
0-framed surgeries on  $g(S^3_1)$ along all $Y_1$ and $Y_2$. 
Then the pass-moves are carried out in all $B^3$. 

Therefore $\Lambda$ is an embedded 3-sphere in $S^5$. 
Furthermore  $\Lambda\cap g(S^3_2)$ in  $\Lambda$ is equal to the 1-link $L_1$   
and $\Lambda\cap g(S^3_2)$ in  $g(S^3_2)$ is equal to the 1-link $L_2$.  

Put $h\vert_{ S^3_1}$ so that  $h(S^3_1)$=$\Lambda$. 
By the above construction, $h$ satisfies the conditions (1)(3)(4)(5) in Lemma 4.8.

We prove: 

{\bf Lemma 4.9} 
  {\it 
$h\vert(S^3_1)$ is equal to the trivial 3-knot.  
}

{\bf Proof.}  
By the construction, $h(S^3_1)$ bounds a 4-manifold 
represented by the disjoint union of some copies of 
the  framed link in 
Figure 9. 

$$\mathrm{Figure 9}$$

See [10] 
for dot-circles. 
Hence this 4-manifold is diffeomorphic to the 4-ball. 
Therefore $h(S^3_1)$ bounds a 4-ball. 
This completes the proof of Lemma 4.9 and therefore Lemma 4.8. 
This completes the proof of Proposition 4.2 and therefore Proposition 4.1.

\section{ The  proof of Theorem 1.3  }

{\bf Lemma 5.1.1 } 
 {\it 
Let $L$ be the trivial 2-component 1-link. 
There exists a self-transverse immersion  
$f:S^3\i B^5$ with the following properties. 

(1)
The singular point set (in $B^5$) is a single circle $C$.  

(2)
$f^{-1}(C)$ in $S^3$ is equivalent to $L$. 

(3)
$f(\overline{S^3-N(L)})$ $\subset \partial B^5$
and 
 $f(\mathrm{Int}N(L))$ $\subset$ $\mathrm{Int}B^5$, where 
 $N({L})$ is a tubular neighborhood of $L$ in $S^3$.  
}

{\bf Lemma 5.1.2 } 
 {\it 
Let $L$ be the Hopf link. 
There exists a self-transverse immersion  
$f:S^3\i B^5$ satisfying with the properties (1), (2) and (3) in Lemma 5.1.1. 
}

{\bf Lemma 5.1.3 }   
  {\it 
Let $L$ be the trivial 1-knot. 
There exists a self-transverse immersion  
$f:S^3\i B^5$ satisfying  with the properties (1), (2) and (3) in Lemma 5.1.1.  }

{\bf Proof of Lemma 5.1.1 }
Take a chart of $(U, \phi)$ of $B^5$ such that  
$\phi (U)$=
$\{(x,y,z,v,t)\vert$ $x,y,z,v \in R,  t\leqq0 \}.$

Put $F$=
$\{(x,y,z,v,t)\vert$ $x,y\in R,  z\geqq0, v=0, t\leqq0 \}.$  

Put $A$=
$\{(x,y,z,v,t)\vert$ $x,y\in R,  z=0, v=0, t\leqq0 \}.$  

We can regard $U$ as the result of rotating $F$ around the axis $A$. 

Take an immersed 2-disc $D$ in $F$ as in 
Figure 10, 11. 


$$\mathrm{Figure 10}$$

$$\mathrm{Figure 11}$$


\vskip3mm
{\bf Note. }    
In $F$ the following hold. 

(1) $\{(x,y,z,v,t)\vert$ $x,y\in R,  z\geqq0, v=0, t=0 \}\cap D$ 
  is diffeomorphic to $\overline{D-(two 2-discs)}$.

(2) $\{(x,y,z,v,t)\vert$ $x,y\in R,  z\geqq0, v=0, t<0 \}\cap D$ 
  is a union of the interior of two 2-discs, 
  where the intersection of the two 2-discs is one point. 
  The point is $p=(1,0,0,0,-1)$. 
  
(3) Let $-1<t'<0$.

 $\{(x,y,z,v,t)\vert$ $x,y\in R,  z\geqq0, v=0, t=t' \}\cap D$ 

is the Hopf link in 

 $\{(x,y,z,v,t)\vert$ $x,y\in R,  z\geqq0, v=0, t=t' \}$. 

(4) $\{(x,y,z,v,t)\vert$ $x,y\in R,  z\geqq0, v=0, t=-1 \}\cap D$ 
  is a union of two circles, 
  where the intersection of the two circles is one point.  
The point is  $p=(1,0,0,0,-1)$. 
 
(5) Let $-2<t'<-1$. 

 $\{(x,y,z,v,t)\vert$ $x,y\in R,  z\geqq0, v=0, t=t' \}\cap D$ 

is the trivial link in 

 $\{(x,y,z,v,t)\vert$ $x,y\in R,  z\geqq0, v=0, t=t' \}$.

(6) $\{(x,y,z,v,t)\vert$ $x,y\in R,  z\geqq0, v=0, t=-2 \}\cap D$ 

is a disjoint union of two 2-discs. 

\vskip3mm


As we rotate $F$ as above, we rotate $D$ as well. 
We obtain an immersed 3-sphere $X$. 

Take $f$ so that $f(S^3)$=$X$. 

This completes the proof of Lemma 5.1.1.

{\bf Proof of Lemma 5.1.2 }
Take $B^5, F, A,$ and $D$ as above. 

Put $G$=
$\{(x,y,z,v,t)\vert$ $x=0,y=0, z\geqq0, v=0, t\leqq0 \}.$

In Figure 10, 11 
we suppose: 
If, in $F$, we rotate $D$ around $G$ by any angle, then $D\cap A$ does not change. 

As we rotate $F$ as above, 
we rotate $D$ around $A$ so that 
in $F$ we rotate $D$ around $G$ one time. 
We obtain an immersed 3-sphere $X$. 

Take $f$ so that $f(S^3)$=$X$.     

This completes the proof of Lemma 5.1.2.

{\bf Proof of Lemma 5.1.3 }
Take $B^5, F, A, D, $ and $G$ as above.

In Figure 10, 11 
we suppose: 
If we rotate $D$ around $G$ half time, 
then the resultant $D$ coincides with the original $D$.

As we rotate $F$ as above, 
we rotate $D$ around $A$ so that 
in $F$ we rotate $D$ around $G$ half time. 
We obtain an immersed 3-sphere $X$. 

Take $f$ so that $f(S^3)$=$X$. 

This completes the proof of the proof of Lemma 5.1.3.

{\bf Lemma 5.2.1 } 
 {\it 
Let $L$ be a 2-component 1-link whose linking number is even. 
There exists a self-transverse immersion 
$g:S^3\i \natural^l S^2\x B^3$, 
for a non-negative integer $l$,  with the following properties. 

(1)
The singular point set (in $\natural^l S^2\x B^3$) is a single circle $C$.  

(2)
$g^{-1}(C)$ in $S^3$ is equivalent to $L$. 

(3)
$g(\overline{S^3-N(L)})$ $\subset \partial(\natural^l S^2\x B^3)=$
$\sharp^l S^2\x S^2$  
and 
 $g(\mathrm{Int}N(L))$ $\subset$ $\mathrm{Int}(\natural^l S^2\x B^3)$, where 
 $N({L})$ is a tubular neighborhood of $L$ in $S^3$.  

Here, we suppose that $\natural^0 S^2\x B^3$ means $B^5$.  
}

{\bf Lemma 5.2.2 }  
 {\it 
Let $L$ be a 2-component 1-link whose linking number is odd. 
There exists a self-transverse immersion 
$g:S^3\i \natural^l S^2\x B^3$, 
for a non-negative integer $l$,  
satisfying with the properties (1), (2) and (3) in Lemma 5.2.1. 
}

{\bf Lemma 5.2.3 } 
 {\it 
Let $L$ be a 1-knot.  
There exists a self-transverse immersion  
\nl$g:S^3\i \natural^l S^2\x B^3$, 
for a non-negative integer $l$,  
satisfying with the properties (1), (2) and (3) in Lemma 5.2.1. 
}

We prove:  

{\bf Claim 5.2.4} 
 {\it 
Lemma 5.2.1, 5.2.2, and 5.2.3 imply Theorem 1.3. 
}

There is an embedding 
$h:\natural^l S^2\x B^3\e S^5$. 
Then $h\circ g$ is an immersion in Theorem 1.3. 
This completes the proof.

In order to prove Lemma 5.2.1, 5.2.2, and 5.2.3, 
we review the $\sharp$-moves.  
See  [14] and [7] for detail.

{\bf Definition }  ([14]) 
Two 1-links are  {\it $\sharp$-move equivalent}
 if one is obtained from the other 
by a sequence of  $\sharp$-moves. 
See 
Figure 12 
for an illustrations of the $\sharp$-move. 
The  $\sharp$-move is different from the pass-move by the orientation.  


$$\mathrm{Figure 12}$$

In each 3-ball in 
Figure 12 
there are four arcs. 
Each of four arcs may belong to different components of the 1-link. 
We do not assume the four arcs belong to one component of the 1-link.

The following proposition 5.3 is proved in the appendix of [14].

{\bf Proposition 5.3 } ([14])
 {\it 
(1)
Let $L$ be a 2-component link. 
Then $L$ is  $\sharp$-move  equivalent to the trivial link 
if and only if the  linking number is even.  

(2)
Let $L$ be a 2-component link. 
Then $L$ is  $\sharp$-move  equivalent to the Hopf link   
if and only if the  linking number is odd.  

(3)
Any 1-knot is  $\sharp$-move equivalent to the trivial knot.
}

We have:

{\bf Proposition 5.4 } ([14])
 {\it 
Let $L_1$ and $L_2$ be 1-links. 
Suppose that $L_1$ is $\sharp$-move equivalent to  $L_2$. 
 Take $L_1$ in $S^3$. 
Then there exists a disjoint union of 3-balls $B^3$ in $S^3$ 
such that $L_1\cap B^3$ is as in 
Figure 12  
and that  the $\sharp$-moves in all $B^3$ change $L_1$ to $L_2$.  
}

Let $L_1$ and $L_2$ be 1-links. 
Suppose $L_1$ is $\sharp$-move equivalent to $L_2$. 
Take $L_1$ in $S^3_1$.  
Take a disjoint union of 3-balls $B^3$ in $S^3$ as in Proposition 5.3.     
Take a 1-link ($Y_{1}$, $Y_{2}$) in each $B^3$ as in 
Figure 13.
Suppose that $Y_i$ bounds a 2-disc $B^2_i$ as in 
Figure 13.   


$$\mathrm{Figure 13}$$

By considering the Kirby moves of framed links (in [10]), 
it is easy to prove:

{\bf Lemma 5.5 }  
 {\it 
We carry out 0-framed surgeries on all $Y_{1}$ and all $Y_{2}$.    
After these surgeries,  we have the following.  
(1)
each $B^3$ becomes a 3-ball again,  
(2)
$S^3$ becomes  a  3-sphere  again. 
(3) $L_1$ in the old sphere $S^3$ changes to $L_2$ in the new 3-sphere. 
}

We go back to the proof of Lemma 5.2.1, 5.2.2, 5.2.3.

{\bf Proof of Lemma 5.2.1 ( resp. 5.2.2, 5.2.3 )} 

By Proposition  5.3, 
a sequence of $\sharp$-moves changes $L$ into 
the trivial 2-component 1-link 
( resp. the Hopf link, the trivial 1-knot ).

As in Lemma  5.1.1 ( resp. 5.1.2, 5.1.3 ),  take $f:S^3\i B^5$. 
See the 1-link $f^{-1}(C)$ in $S^3$. 
It is equivalent to the 1-link $L$.  
As in Proposition 5.4, take 3-balls $B^3$. 
As in Lemma 5.5, take $Y_1$ and $Y_2$ in $B^3$.  
Suppose that $f(Y_1)$ and $f(Y_2)$ are in $\partial B^5$.

Attach 4-dimensional 2-handles $h^2$ to 
$f(S^3-N(L))$  in $\partial B^5$ along $f(Y_i)$ with 0-framing.  
Here, 0-framing means the following. 
When attaching the 4-dimensional 2-handles $h^2$ to $f(S^3-N(L))$,  
we can attach 4-dimensinal 2-handles to $S^3$ naturally. 
These attaching maps are 0-framing. 

When attaching the 4-dimensional 2-handles $h^2$ to $f(S^3-N(L))$,  
we can attach 5-dimensional 2-handles  $h^2\x[-1,1]$ 
to $B^5$ along $f(Y_i)$ naturally.  
Of course the attached parts are in  $\partial B^5$. 
We obtained a 5-manifold. Call it $M$.

Put $P=$ $f(S^3)-(f(S^3)\cap\partial h^2)$. 
Put $Q=$ $\partial h^2-(f(S^3)\cap\partial h^2)$. 
Put $\Sigma=P\cup Q$. 
Then $\Sigma$ is an immersed 3-sphere in $M$. 

We prove Lemma 5.6. 
Before the proof of Lemma 5.6, we prove Lemma 5.7.

{\bf Lemma  5.6 }  
 {\it 
Let $M$ be as above. 
$M$=$\natural^l S^2\x B^3$, for a non-negative integer $l$. 
}

{\bf Lemma  5.7 }  
 {\it 
Lemma 5.6 implies Lemma 5.2.1, 5.2.2 and 5.2.3.  
}

{\bf Proof of Lemma 5.7.  }
Take a self-transverse immersion  $g:S^3\i M=\natural^l S^2\x B^3$
so that  $g(S^3)=\Sigma$.

{\bf Proof of Lemma 5.6.}
Recall:
$J=$ (one 5-dimensional 0-handle)$\cup$(one 5-dimensional 2-handle) is 
$S^2\times B^3$ or $S^2\tilde\times B^3$. 
 $J$ is $S^2\times B^3$ 
if and only if the attaching map of the 2-handle is spin preserving diffeomorphism map. 

It suffices to prove that 
the attaching maps of the 5-dimensional 2-handles are 
spin-preserving diffeomorphism maps. 

We give a spin structure on $f(S^3-N(L))$ 
from the unique one on $\partial B^5$.  

We give a spin structure on $f(B^2_i)\cap\partial B^5$ 
from the spin structure on $f(S^3-N(L))$.  

We give a spin structure $\xi$ on $f(Y_i)$
 from the spin structure on $f(B^2_i)\cap\partial B^5$.  

Put $S^1_a\amalg S^1_b$= $f(B^2_i)\cap\partial N(L)$.  
We give a spin structure $\alpha$ (resp. $\beta$)  
on $S^1_a$ (resp.  $S^1_b$ ) 
 from the spin structure 
on $f(B^2_i)\cap\partial B^5$.  
By the construction,  $\alpha$ and  $\beta$  are the $S^1_{Lie}$ spin structure.

Since $(f(Y_i), \xi)$ is spin cobordant to 

$(S^1_a,$ the $S^1_{Lie}$ spin structure $)$ $\amalg$ 
$(S^1_b,$ the $S^1_{Lie}$ spin structure $)$,   

$\xi$ is  the $S^1_{bd}$ spin structure.

See [10] 
for the $S^1_{bd}$ spin structure and the $S^1_{Lie}$ spin structure.  

This completes the proof of Lemma 5.6. 

This completes the proof of Lemma 5.2.1 ( resp. 5.2.2, 5.2.3 ).   
By Claim 5.2.4, Theorem 1.3 holds.

\np

\noindent{\large{\bf Appendix. The proof of  Proposition 2.1 and 2.2 }}

Firstly we prove Proposition 2.1. 

Let $L=(K_1,...,K_m)$ be a 1-link. 
Let $F$ be a Seifert surface. 
Let $K_i\x[0,1] (\subset F)$ be a collar neighborhood of $K_i$ in $F$. 
Let  $K_i\x\{0\}$=$K_i$. 
Then lk$(K_i, L-K_i)$=lk$(K_i\x\{0\}, K_i\x\{1\})$. 

See [10] 
for spin structures.   

Let $\varepsilon^2$ be the trivial bundle over $S^1$. 
Let $(e^1_p, e^2_p)$ (for each $p\in S^1$) denote the trivialization. 
Let $e^0_p$ (for each $p\in S^1$) denote the trivialization on $TS^1$.  
Spin structures on $S^1$ are defined by spin structures on 
$TS^1\oplus\varepsilon^2$. 
[10]
defines that the spin structure defined by $(e^0_p, e^1_p, e^2_p)$
is the {\it $S^1_{Lie}$ spin structure}. 
The other is the {\it $S^1_{bd}$ spin structure}.

Let $(M, \sigma)$ and $(N,\tau)$ be spin manifolds. 
Let $f:M\rightarrow N$ be an orientation preserving diffeomorphism. 
$f$ is called a {\it spin preserving diffeomorphism}  
if  
$df\oplus id:TM\oplus\varepsilon^p  \rightarrow TN\oplus\varepsilon^p $
carries 
$\sigma$ to $\tau$, 
$id$ carries the trivialization to the trivialization.

We give $F$ a spin structure $\alpha$ induced from 
the unique spin structure $\beta$ on $S^3$. 
We give $K_i$ a spin structure $\gamma$ induced from $\alpha$.  
Let $K_i\x D^2$ be the tubular neighborhood of $K_i$ in $S^3$. 
Regard  $K_i\x D^2$ as the product $D^2$ bundle over $K_i$. 
We give $K_i\x D^2$ a trivialization such that 
$e^1_p$ is in  $K_i\x[0,1]$ and that 
$e^2_p$ is perpendicular to $K_i\x[0,1]$ for each $p\in K_i$.
Then we can regard $\alpha$ as a spin structure on 
the trivialized bundle $TK_i\oplus (K_i\x[0,1])$.

Let $G$ be any Seifert surface of $K_i$. 
Let $g$ be a spin structure on $G$ 
induced from the unique one $S^3$. 
Let $\mu$ be a spin structure on $K_i$ induced from $(G,g)$. 
Since $\partial(G, g)=(K_i, \mu)$, $[(K_i, \mu)]=0$.  
We give $K_i\x D^2$ a trivialization such that 
$g^1_p$ is in $G$ and $g^2_p$ is perpendicular to $G$. 
Then the spin structure on $K_i$ induced by 
the framing $(e^0_p, g^1_p, g^2_p)$ is $\mu$.

The homotopy class of the framing $(e^0_p, e^1_p, e^2_p)$ coincides with 
the homotopy class of the framing $(e^0_p, g^1_p, g^2_p)$ 
if and only if 
lk$(K_i\x\{0\}, K_i\x\{1\})$ is even.

Therefore Proposition 2.1 holds.

Next we prove Proposition 2.2. 

Under the above condition, 
suppose $L$ is a proper 1-link. 

Let $\hat F$ be the closed oriented surface $F\cup(\cup^{m}_{i=1}D^2)$.
Since $\gamma$ is the $S^1_{Lie}$ spin structure, 
$\alpha$ extends to $\hat F$, say $\hat\alpha$.  
Let $a_1,...,a_g,b_1,...,b_g$ be circles in $F$ representing 
symplectic basis of 
$H_1(\hat F; Z)$. 
Let $a_i\x[-1,1]$, $b_i\x[-1,1]$ be the tubular neighborhood in $F$. 
Put mod 2 lk ($a_i\x\{-1\}, a_i\x\{1\}$)=$x_i$  
and mod 2 lk ($b_i\x\{-1\}, b_i\x\{1\}$)=$y_i$. 
Recall Arf$L=$$\Sigma x_i\cdot y_i$.

We give a spin structure $\sigma_i$, (resp. $\tau_i$) 
on $a_i$ (resp. $b_i$). 
Then [($a_i, \sigma_i$)]$\in\Omega_1^{spin}=$ $x_i$ and 
 [($b_i, \tau_i$)]$\in\Omega_1^{spin}=$ $y_i$. 
By P.36 of [10],
[($\hat F, \hat\alpha$)]$\in\Omega_2^{spin}$ is $\Sigma x_i\cdot y_i$ .

Hence Proposition 2.2 holds.

\np
\footnotesize{

\newpage

\unitlength 0.1in
%

The shaded part is $F\cap B^3$.

\newpage

\unitlength 0.1in
%
%

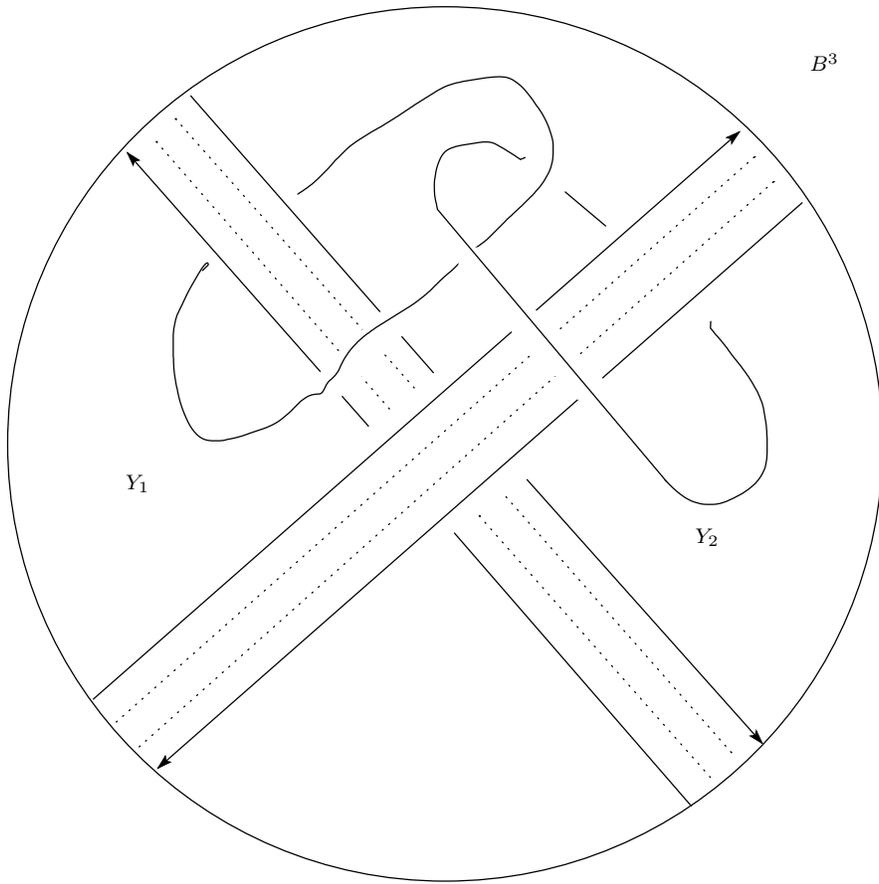
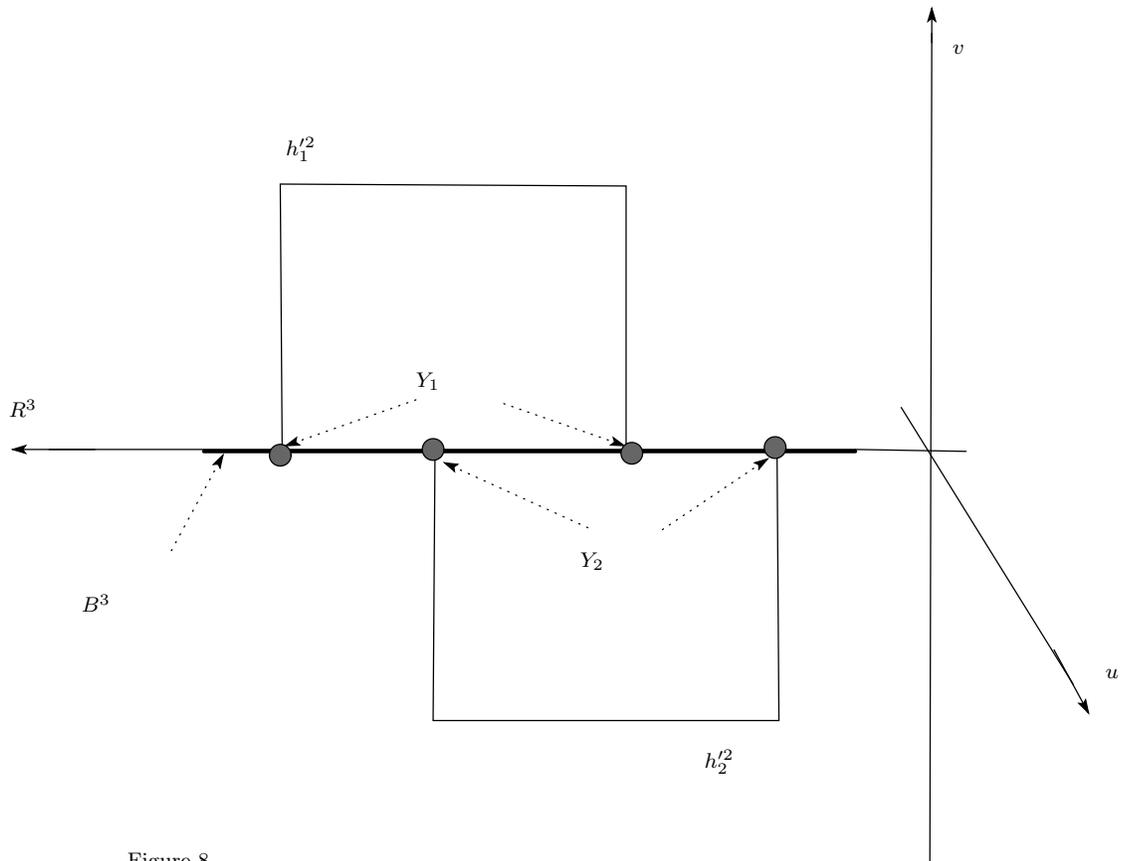
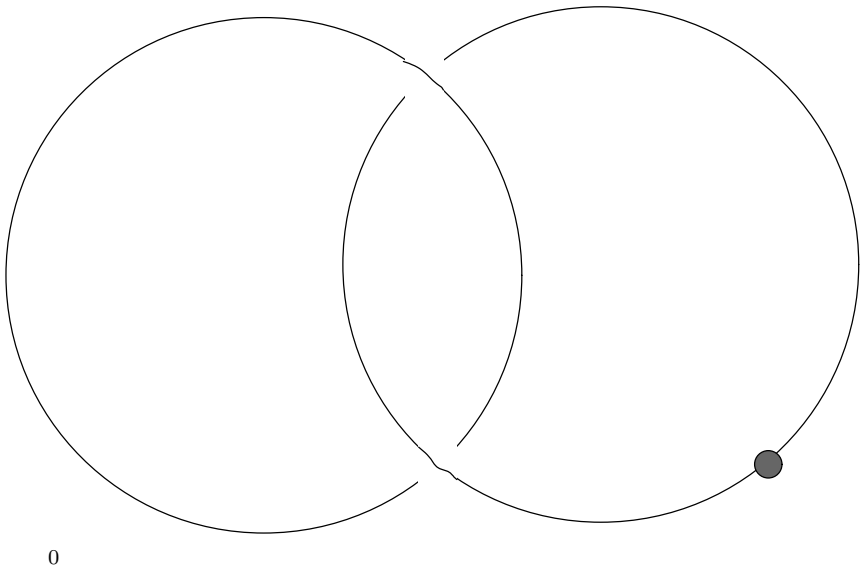
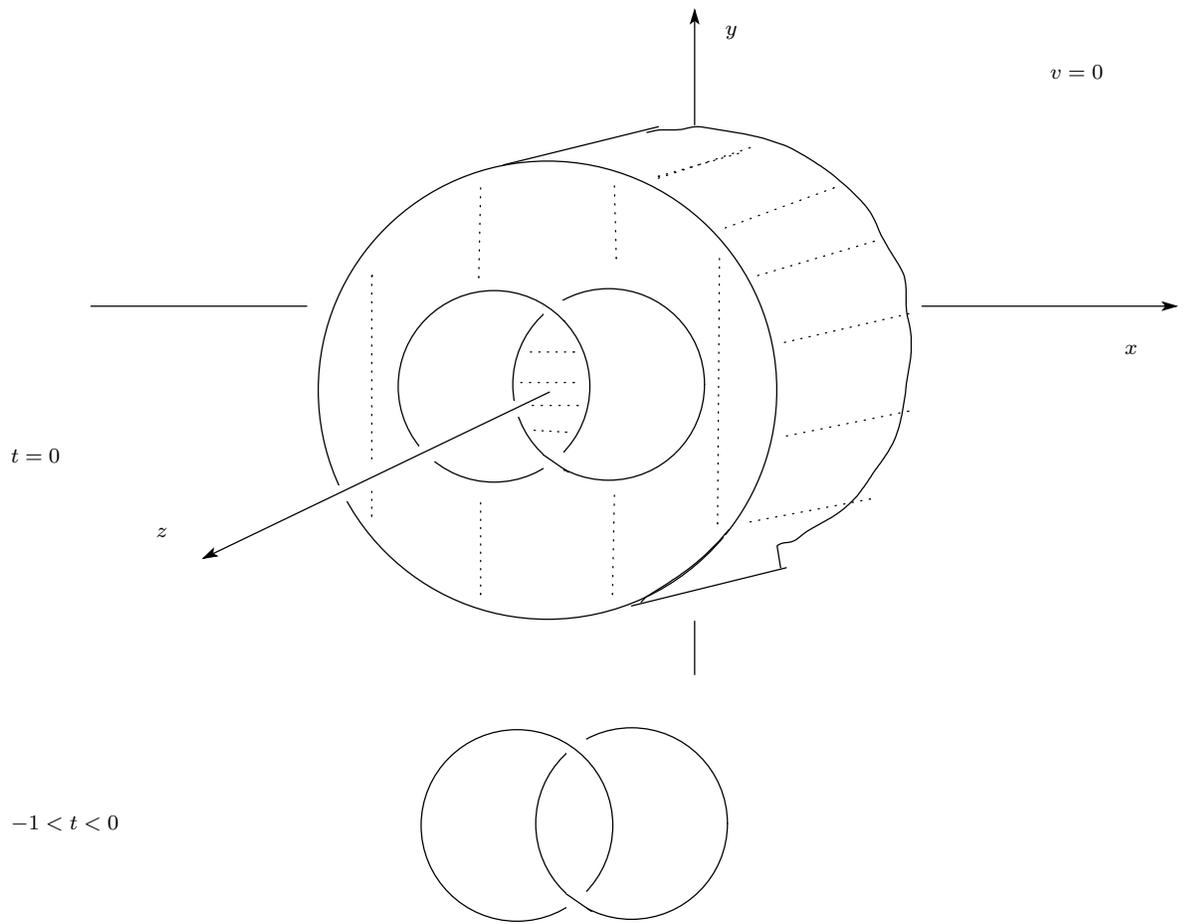
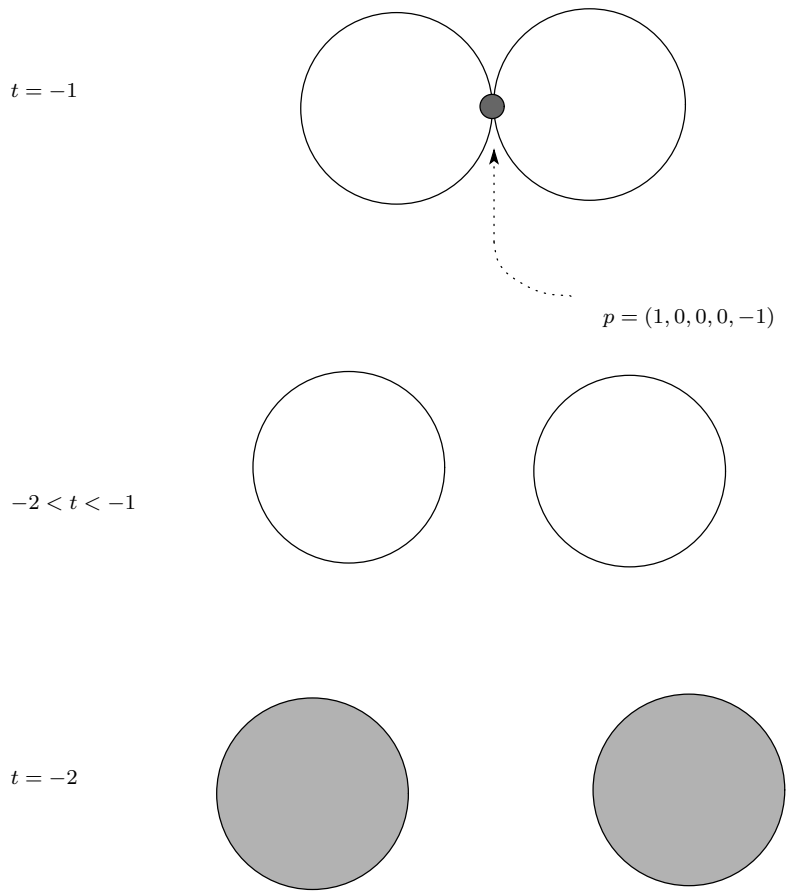
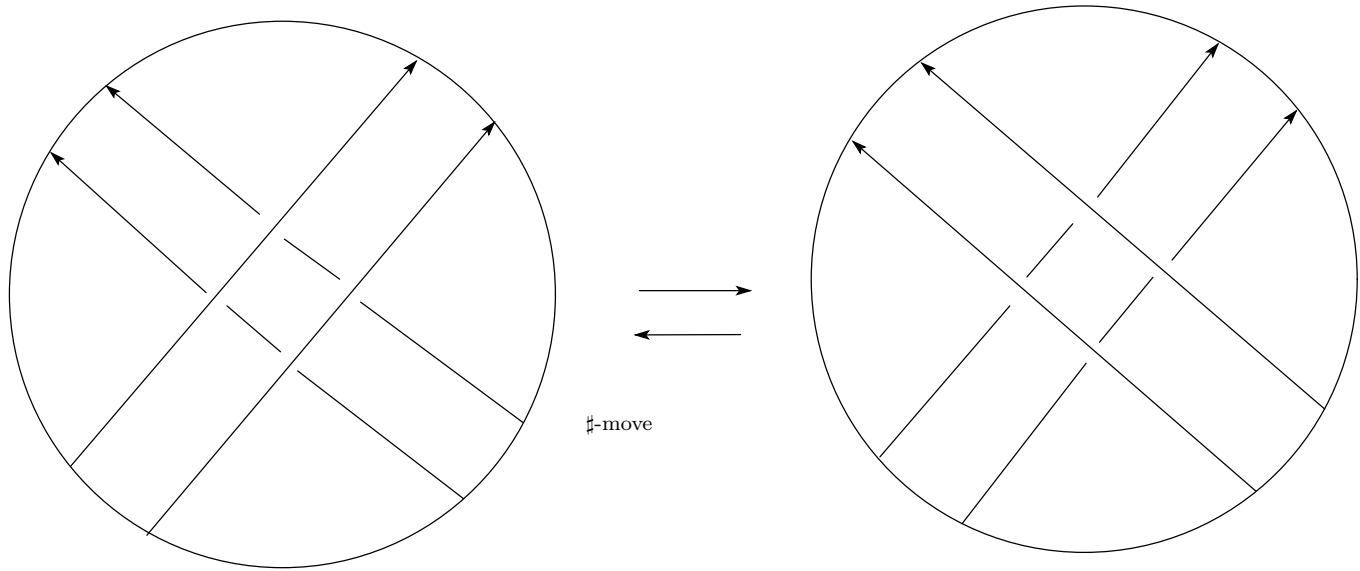
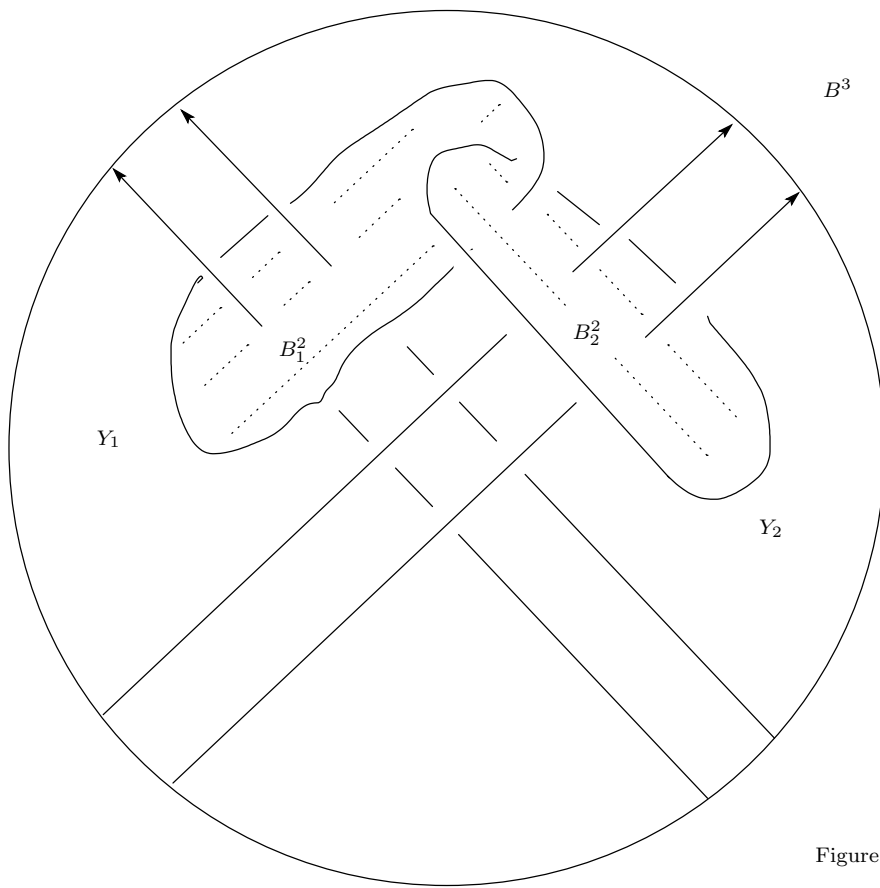
\end{document}